\documentclass[number,citesort,seceqn,dvips]{arxbj}
\usepackage{graphicx}

\aid{0}
\volume{17}
\issue{3}
\pubyear{2011}
\firstpage{895}
\lastpage{915}
\doi{10.3150/10-BEJ303}

\makeatletter
\newcommand{\eqref}[1]{(\ref{#1})}
\def\F{\mathcal{F}}
\def\G{\mathcal{G}}
\def\m{^{(m)}}
\newtheorem{theorem}{Theorem}
\newproclaim{definition}{Definition}
\newremark{remark}{Remark}
\makeatother

\begin{document}
\begin{frontmatter}

\title{A martingale approach to continuous-time marginal structural models}
\runtitle{Martingales and marginal structural models}

\begin{aug}
\author{\fnms{Kjetil} \snm{R{\O}ysland}\corref{}\ead[label=e1]{kjetil.roysland@medisin.uio.no}}
\runauthor{K. R{\o}ysland}
\address{Department of Biostatistics,
Institute of Basic Medical Sciences,
University of Oslo,
Boks 1122 Blindern,
0317 Oslo,
Norway. \printead{e1}}
\end{aug}

\received{\smonth{3} \syear{2009}}
\revised{\smonth{6} \syear{2010}}

%
\begin{abstract}
Marginal structural models were introduced in order to provide
estimates of causal effects from interventions based on observational
studies in epidemiological research. The key point is that this can be
understood in terms of Girsanov's change of measure. This offers a
mathematical interpretation of
marginal structural models that has not been available before.
We consider both a model of an observational study and a model of a
hypothetical randomized trial. These models correspond to different
martingale measures -- the observational measure and the randomized
trial measure -- on some underlying space. We describe situations where
the randomized trial measure is absolutely continuous with respect to
the observational measure. The resulting continuous-time likelihood
ratio process with respect to these two probability measures
corresponds to the weights in discrete-time marginal structural models.
In order to do inference for the hypothetical randomized trial, we can
simulate samples using observational data weighted by this likelihood ratio.

\end{abstract}

%

\end{frontmatter}

\section{Introduction}
We will consider the following scenario: A patient has a disease. In
order to avoid an event (e.g., death), a specific treatment can
be given. The given treatment will typically depend on the patient's
previous health condition.

We would like to estimate the effect of a given treatment on the time
until the occurrence of the event.
A natural way to do so is to implement some sort of randomized trial.
This means that we would have to set up an experiment on a group of
patients where the treatment was initiated
by randomization independently of each patient's previous health condition.
Such a study typically requires significant resources that may not be
available.
In order to take advantage of another type of data, we could try to
base our estimates of the treatment effect on an observational
study. Suppose we have observations of a group of patients where the
given treatments were chosen by doctors.
As a first attempt, one could try to compute the relative short-term
risk between the group given treatment and the group not given
treatment at a particular time. This could be done using Cox
proportional hazards regression techniques.
However, such a naive analysis would most likely introduce a bias
compared to the estimate based on the randomized trial.
The reason is that the health condition of the patient not already in
treatment will be a predictor of both treatment and death; that is, it
is likely to be a confounder \cite{Sterne}.

We can easily imagine two opposite scenarios where this confounder
would complicate estimates:
Due to considerable costs, reduced life quality or possibly drug
resistance, one could decide that
the treatment should not be initiated until the patients are
sufficiently ill.
A naive\vadjust{\goodbreak} marginal analysis based on data from an
observational study would then quickly lead us to believe that
the treatment effect was less than the true treatment effect.
Conversely, if we decided only to initiate treatment for patients with
good health conditions and not for the ones with poor conditions, then
a naive marginal analysis would quickly lead us to believe that the
treatment effect was better than the true treatment effect.

In order to solve this problem, one might suggest that we compute an
estimate of the treatment effect conditionally
on the health condition of the patient. However, in several situations,
it is likely that the previous treatment will improve the patient's
general intermediate health condition. This improvement will in itself typically
postpone the time of death. The conditional effect estimate we
described would only incorporate the direct treatment effect, not the effect
that is due to an improvement of the patient's intermediate health condition.

There is also another source of bias that we have to consider in order
to lay hands on the causal effect of treatment; that is, censoring.
We assume that a patient may drop out of the study at a time and not
return; that is, we have right censoring.
The given treatment, calender time and the patient's health condition
might lead to such a drop out. If we do a naive analysis based on the
patients that are still in the study, then we introduce a selection
bias \cite{hernan2}.

We are forced to move outside the standard Cox regression framework
since we have to deal with the mentioned
time-dependent confounder effects due to a patient's underlying health
condition.
In order to provide a meaningful estimate of the treatment effect, with
a simple interpretation,
we could try to construct a rich model that also describes the dynamics
of the underlying biological processes.
Such mechanisms are likely to be very complicated and there might not
be sufficient knowledge or data available.
For this reason we could
try to fit a marginal model of a suitable randomized trial for our
scenario. This will be our strategy in what follows.

One attempt to provide a marginal estimate of the causal treatment
effect this way is presented by Robins in
\cite{Robins1}. This method uses marginal structural models and
relies on the additional assumption that there are no unmeasured
confounders; that is, there does not exist an unobserved process that
is a
predictor of both censoring and treatment, both censoring and the event or
both treatment and the event, given the observed covariates.
If every such process is measured, then the marginal structural model (MSM) approach provides a
proper adjustment
of the marginal effect estimates.
The idea is to apply some clever weights to the observations. This
weighting results in a pseudo-population that is different from the
observed population. The key property of this pseudo-population is that
the selection bias and the treatment confounding due to the patient's
health condition become negligible.
Now, one can proceed with a weighted Cox regression to obtain a marginal
estimate of the effect of treatment. The method has been used several
times on epidemiological studies. In \cite{Zidovudine} it was used to
estimate the effect of Zidovudine on the survival of HIV-positive
men in the Multicenter AIDS Cohort Study. Moreover, the method was also
used in \cite{Sterne} to give an estimate of the hazard ratio for the
effect of highly active antiviral treatment (HAART) on progression to
AIDS or death for HIV patients in Switzerland.

The method introduced by Robins deals with longitudinal data in
discrete time.
We will consider continuous-time versions of the marginal structural
models for event history data. The idea is to characterize
reasonable models of a randomized trial, the randomized trial
measures, using martingale theory.
This offers a mathematical interpretation of marginal structural models
that has not
been available before.

We characterize a class of of reasonable models
of randomized trials in terms of local independence.
Such a model corresponds to a particular martingale measure.
The continuous-time likelihood ratio process between this measure and
the observational
probability measure corresponds to the weights in
a discrete-time marginal structural model.
In order to do inference for this new measure, we can simulate
samples using the observed data weighted by this likelihood ratio.

Another approach to causal inference within our scenario is to use the
so-called structural nested models. These models were also introduced
by Robins; see \cite{Rob3,Rob4}. Lok has developed
continuous-time versions of such models using counting processes and
martingale theory; see~\cite{Lok}.

\section{Observable processes and local independence}
Before we come to the main results, we will spend some time
establishing terminology. Even if the mathematics involved is fairly
standard stochastic process theory, it is perhaps not so commonly used
in event history analysis. A very good background reference on
stochastic processes that we will use frequently is \cite{JacodShiryaev}.

\subsection{Observable processes}
In Section \ref{patientmodel} we will consider a stochastic model of a single
patient. There are typically many factors that are
important for describing how the disease of that individual
develops in time.
We will consider models where all the possible observations of one
patient are represented by stochastic integrals against Poisson
processes. More formally, let $d, n \in\mathbb N$ and consider a
probability space
$(\Omega, \F, Q)$ with mutually orthogonal counting processes $N_t^1,
\ldots, N^n_t$ on the interval
$[0, T]$ and a filtration
$\{\F_t\}_t$ that is generated by
their joint history and some initial information~$\F_0$.
The counting processes are assumed to be Poisson processes in the sense that
\[
\overline{N}{}^1_t := N_t^1 - t, \ldots,\overline{N}{}^n_t := N_t^n - t
\]
define $Q$-martingales, the compensated Poisson processes.
The probability measure $Q$ will only play a role as a reference
measure, as
we will mainly be interested in probability measures that are
absolutely continuous with respect to $Q$. This will sometimes be
referred to as the
\textit{Poisson measure}.

We let
$H$ be a bounded and $\F_t$-predictable $d \times n$-matrix-valued
process and
let $X_0$ denote a bounded $\F_0$-measurable random vector. Now,
define the $d$-dimensional \textit{observable process}:
\[
\vspace*{-2pt}
X_t := X_0 + \int_0 ^t H_s \,\mathrm{d}N_s.
\vspace*{-2pt}
\]
All the possible observations of a patient in our approach will be processes
of this form. Counting processes are trivially included in this class,
but we also allow slightly more complicated jump processes. One
example could be measurements of blood values. Each time the blood
value is updated, it would be given by a jump and correspond to a jump
time of the underlying counting process. The size and direction of
the jump would then be given by the value of the predictable integrand
$H$ at the jump
time.

\subsection{Separability}
We will say that two observable processes $X$ and $Y$ are \textit{separable}
if
they allow the representations:
\begin{eqnarray*}
\vspace*{-2pt}
X_t & =& X_0 + \int_0 ^t H^X_s \,\mathrm{d}N_s^X, \\
Y_t & =& Y_0 + \int_0 ^t H^Y_s \,\mathrm{d}N_s^Y,
\vspace*{-2pt}
\end{eqnarray*}
where $N^X$ and $N^Y$ are independent components of the multivariate
process $N$,
$X_0$ and $Y_0$ are bounded $\F_0$-measurable random vectors and $H^X$
and $H^Y$ are bounded matrix-valued processes
that are predictable with respect to the histories of $N^X$ and $N^Y$,
respectively.
Separability is a technical assumption that provides well-behaved
factorizations of
likelihoods. This is used in the proof of Theorem \ref{dsepcens}.
Heuristically, it means that the processes $X$ and $Y$ do reflect
different random phenomena. Separability is even stronger than orthogonality
since the processes are independent with respect to the Poisson
measure $Q$. However, since we will deal with other probability measures
that are absolutely continuous with respect to the Poisson measure,
separable processes can not necessarily be treated as
independent.

\subsection{A martingale measure}
As we mentioned earlier, our samples will consist of paths of
observable processes.
These samples will be distributed according to
some probability measure $P$
such
that a given family of predictable and non-negative processes
define the jump intensities for
$N^1, \ldots, N^n$ with respect to $\F_t$.
Since we assume that the observations are distributed according to
$P$, we will refer to such a measure as an \textit{observational measure}.

More
formally,
we let $\lambda^1, \ldots, \lambda^n$ be non-negative $\F_t$-predictable
processes
and we
assume that $P$ is a probability measure such that:
\begin{enumerate}[(1)]
\item[(1)]
$P$ and $Q$ coincide on $\F_0$,
\item[(2)]
$P \ll Q$, i.e., $P$ is absolutely continuous with respect to the
Poisson measure,
\item[(3)] The equation
\[
M^i_t := N^i_t - \int_0 ^t \lambda_s^i \,\mathrm{d}s
\]
defines a square-integrable $P$-martingale with respect to $\mathcal
F_t$ for every $i$.
\end{enumerate}
These properties characterize the probability measure $P$ uniquely if
such a measure exists, \cite{JacodShiryaev}, Theorem III 1.26.

\subsection{Non-influence}
We will need a notion of non-influence between observable
processes. There are several formal definitions that are meant to
capture this; see \cite{Florens}.
Independence, or even
conditional independence, is too strong to be of interest for the
method we have in mind. The non-influence relation we will
consider is \textit{local independence}.
Heuristically, a process $X$ is locally independent of a process $Y$
if information about the past
of $Y$ does not contribute to a better prediction of the short-term
behavior of $X$.

In the setting of event history analysis, this concept has been studied
thoroughly by Didelez~\cite{Didelez}. Schweder \cite{Schweder}
used this concept in a study of composable Markov processes.
Aalen \textit{et al.} made use of local independence in order to
study the effect of menopause on the risk of developing a
certain skin disease in \cite{aalen1980}.

\subsection{Local independence}

Let $X, Y, Z$ be observable processes that are mutually separable.
The processes $X_t - X_0,\break Y_t- Y_0$ and $Z_t - Z_0$ are obviously
independent with respect to the
probability measure~$Q$. However, the
situation is typically more complex with respect to the measure~$P$, since the jump intensities
$\lambda^1_t, \ldots, \lambda^n_t$ could depend on all the information
in $\F_{t-}$. We therefore introduce the following concept.
\begin{definition}
Let $\F_t^{X, Y, Z}$ denote the filtration generated by $N^X, N^Y,
N^Z$ and let
$\F_t^{X, Z}$ denote the filtration generated by $N^X$ and $N^Z$.
We say that $X$ is \textit{locally independent} of~$Y$, given $Z$,
if there exists an $\F_t^{X, Z}$-predictable process $\mu$
such that
\[
N^X_t - \int_0 ^t \mu_s \,\mathrm{d}s
\]
defines a
local $P$-martingale with respect to $\F_t^{X, Y, Z}$. If this is the case,
then we write:
\[
Y \nrightarrow X | Z.
\]
\end{definition}
%
\subsection{Independent censoring}
Local independence generalizes a much-used concept in event history
analysis, that is, \textit{independent censoring}.
Suppose we can follow a group of individuals in a clinical trial.
We would like to compute the probability for an individual to survive
longer than time
$t$. However, an individual might be censored at some time before the
event due to the end of
the study or a ``drop-out''.
Inference is much simpler if the censoring does not influence
the instantaneous risk of the event.
Therefore, it is common to assume independent censoring.
This means that an individual at risk has the same
instantaneous risk of an event as he would in
the situation without censoring. More formally, this means that if
$T_D$ is the time of the event and $T_C$ is the time of censoring, then
the compensator of the process $D_t := I( t \geq T_D)$ with respect
to the joint event and censoring history only depends on the event
history.
This is essentially the same as saying that $D$ is locally independent
of the process defined by $C_t := I( t \geq T_C) $.

\subsection{Local independence before a stopping time}

Sometimes we may not be interested in dependencies that are considered
trivial. This could be dependencies due to
an absorbing state (e.g., death).
We will see that we can rule out such trivial dependencies if
we consider local independence before a stopping time $\tau$.
\begin{definition}
Let $\tau$ be an $\F_t$-adapted stopping time.\vspace*{1pt}
We say that $X$ is locally independent of $Y$ before $\tau$ and given $Z$
if there exists an $\{\F_t^{X, Z}\}_t$-predictable process $\mu$ such that
\[
N^X_{t\wedge\tau} - \int_0 ^{t\wedge\tau} \mu_s \,\mathrm{d}s
\]
defines a local $P$-martingale with respect to $\{\F_t^{X, Y, Z}\}_t$.
If this is the case,
then we write:
\[
Y \nrightarrow_\tau X | Z.
\]
\end{definition}

If we let
$\tau$ denote the time of the first jump of $N^Y$
then it is not very hard to see, using the explicit representation of
$\F_t^{X, Y, Z}$-predictable processes in
\cite{Bremaud}, Theorem A.2, that for every $\F_t$-predictable
process $\gamma$ there exists an $\F_t^{X, Z}$-predictable process
$\tilde\gamma$ such that $\gamma_S \cdot I ( S \leq\tau)= \tilde
\gamma_S \cdot I ( S \leq\tau)$ $P$-a.s. for every $\F_t$-adapted stopping
time $S$. This means that
$
Y \nrightarrow_\tau X | Z
$; that is, stopping at the first jump of $N^Y$ rules out every local dependence
of $Y$.
\subsection{Local independence graphs}
Didelez also considered graphical models based on local independence;
see \cite{Didelez}.
These graphs will prove to be
very useful in order to represent complex models.

\begin{definition}
We say that a directed graph $G = (E, V)$ is a local independence
graph if the vertexes correspond to observable processes that are
mutually separable and
such that
\[
(X, Y) \notin E \quad \Longrightarrow\quad  X \nrightarrow_\tau Y | V \setminus\{X,Y\}.
\]
\end{definition}

Several examples of such graphs will appear below.

\section{Models of clinical trials} \label{patientmodel}

\subsection{A patient model}

We will now describe a model of a single patient that participates in a
clinical study.
We suppose that $N$ is of the form $( N^A, N^C, N^D, N^L )'$,
where $N^A, N^C, N^D$ are univariate counting processes and
$N^L$ is a multivariate counting process. These counting processes
count various events that are important for the development of the
disease.

\subsubsection{The event process}
We let $T_D = \inf\{ t > 0 \mid N_t^D = 1 \}$ and let
\[
D_t := \int_0 ^t I( s \leq T_D) \,\mathrm{d}N_s^D
\]
be the \textit{event process}. It jumps from $0$ to $1$ at the
time the event occurs. The event could be death or the progression to
AIDS for an HIV patient.

\subsubsection{Measurements of the underlying biological process}
The state of an underlying biological process reflecting the
patient's health condition at time $t$ is given by
\[
L_t := L_0 + \int_0 ^t H_s ^L \,\mathrm{d}N_s^L,
\]
where $L_0$ is a bounded $\F_0$-measurable random vector and
$H^L$ is a matrix-valued, bounded and $\F^L_t$-predictable
process.
The process $L$ could be measurements of various blood values.

\subsubsection{Right censoring}
We assume that the patient can be right censored,
that is, we will not be able to observe the patient after some stopping
time $T_C$. This can happen because the study ends,
but it can also be a ``drop-out'' due to poor health or recovery.
We assume that $T_C := \inf\{ s > 0 \mid N^C_s \neq0 \}$
and define the censoring process
\[
C_t := \int_0 ^t I( s \leq T_C) \,\mathrm{d}N_s^C.
\]

\subsubsection{The treatment process}
One can switch between two treatments of the patient
at the stopping time $T_A$. This could typically be to initiate
treatment for a patient at risk.
We let $T_A := \inf\{ s > 0 \mid N^A_s \neq0 \}$
and define the treatment process
\[
A_t := \int_0 ^t I( s \leq T_A ) \,\mathrm{d}N_s^A.
\]
This means that the patient will not initially be in
treatment. This somewhat limiting assumption can be dropped, but then
the considerations around the hypothetical randomized trial at
baseline will be much more involved.

\subsection{Local independences with respect to the observational measure}

The process $D$ influences the other processes. However, we consider
these dependencies as trivial.
We will consider local independence before $T_D$, because
then we will automatically have that
$D \nrightarrow_{T_D} A | C \cup L$, $D \nrightarrow_{T_D} C |
A\cup L$ and $D \nrightarrow_{T_D} L | A\cup C$.

We also assume that the censoring does not carry any information
about the short term behavior of the other process
that we would not obtain if we left $C$ out of the analysis. In terms
of local independence this means that
$C \nrightarrow_{T_D} A | D\cup L$, $C \nrightarrow_{T_D} D | A\cup
L$ and
$C \nrightarrow_{T_D} L |A\cup D$.
We summarize these local independencies in the following local
independence graph:
\begin{equation}\label{Graphicalmodel}
\begin{tabular}{c}
$
\includegraphics{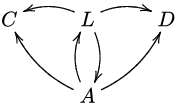}
$
\end{tabular}
\end{equation}

\subsection{Randomized trial measures}
Our ultimate goal is to provide estimates of the causal effect of a particular
treatment based on observations of patients in an observational study.

Hypothetically, one could carry out some randomized trial
where the given treatment did not depend on the previous health
condition of the patient.
If we had observations from such a trial, we could easily provide
simple estimators for the causal treatment effect that
would not require information about the underlying biological
mechanisms. This is, however, not the case for us. So, based on
observations from the observational study, we will try to simulate a
counterfactual or hypothetical randomized trial.
We assume that we have measurements of all the relevant processes and variables.
Especially, we assume that the process $L$ is complete in the sense
that it gives rise to every event
that affects both the short-term behavior of the treatment and the event,
both the censoring and the event or both the censoring and the
treatment, given
the full covariate history. This assumption means that all the confounder
processes are measured and is usually referred to as \textit{no
unmeasured confounders}.

In order to provide causal interpretation of simple estimators,
we should at least require the hypothetical trial to satisfy the
following:
\begin{enumerate}[(1)]
\item[(1)]
Both the underlying biological process and the event process should
dynamically
behave in the same way in the counterfactual trial and the
observational study, given the
full covariate history.
\item[(2)] One should not allow drop-out due to poor health or recovery,
that is, the
censoring should not be directly affected by the
underlying health process, given the event and treatment history.
\item[(3)]
Since we consider time-dependent treatments, we have to generalize the
notion of a randomized trial slightly.
In our counterfactual trial, the patient's previous health condition
or censoring
should not be relevant for the short-term behavior of the
treatment process. Heuristically, this means that the randomization
should act locally in
time.
\end{enumerate}

The counterfactual trial corresponds to a probability measure
$\tilde P$ on the space $( \Omega, \F)$. We will refer to such a
measure as a \textit{randomized trial measure}.
It carries the frequencies of the potential observations in the
counterfactual randomized trial.
The above requirements can now be translated into the following:
\begin{enumerate}[(1)]
\item[(1)] The process $M^i$ is a local $\tilde
P$-martingale with respect to the filtration $\{\F_t\}_t$ for every
$i \in D \cup L$. This means that both the processes
$N^L$ and $N^D$ have the same intensity with respect to the
randomized trial measure $\tilde
P$ as with respect to the observational measure~$P$.
Moreover, we assume that the observational measure and the
randomized trial measure coincide at baseline, that is,
\[
E_P [ H ] = E_{\tilde P} [ H ]
\]
for every
bounded $\F_0$-measurable random variable $H$.

\item[(2)]
The censoring should be locally independent of the
underlying health process $L$ with respect to $\tilde P$, given the
event and treatment history.
\item[(3)]
The treatment process should be locally independent of the underlying
health process $L$ and censoring $C$ with respect to $\tilde P$,
given the event history.
\end{enumerate}

We summarize the local independence structure with respect to the
randomized trial
measure~$\tilde P$ in the following local independence graph:
\[
\begin{tabular}{c}

\includegraphics{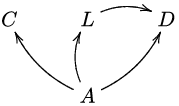}

\end{tabular}
\]

A construction of reasonable randomized trials is given in Theorem
\ref{maintheorem}. Before we come to that, we will consider censoring
in the counterfactual trial.
In order to estimate the total treatment effect, we will consider a
marginal model where $L$ is unobserved.
One natural choice of effect measure in the hypothetical experiment
could be the hazard of
the event process with respect to the filtration $\{\F_t^{A,D}\}_t$.
In order to estimate the hazard with respect to this filtration,
we could try to estimate the hazard of the event before censoring with
respect to
the filtration $\{\F_t^{A,C, D}\}_t$. If, in addition, the event
process was locally independent of the censoring, given
the treatment process, then these hazards would coincide before
censoring; that is, we would have independent censoring. This would
imply that the censoring would not cause bias in the sense that
if we did not pay
attention to the underlying biological process, then
the hazard of the
event would not depend on whether the patient had been censored or
not.
We will see in the next theorem that this is the case for
the randomized trial measures.

\begin{theorem} \label{dsepcens}
If $P$ is a randomized trial measure then we have
$C \nrightarrow_{T_D} A\cup D $, that is, we have independent
censoring in the marginalized model
without $L$.
This gives the following local independence graph:
\[
\begin{tabular}{c}

\includegraphics{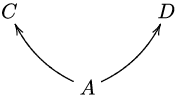}

\end{tabular}
\]
\end{theorem}

\begin{pf}
The likelihood ratio process
\[
S_t := \frac{\mathrm{d} P |_{\F_{t}} }{\mathrm{d} Q
|_{\F_{t}}}
\]
is a $Q$-martingale with respect to $\F_t$. This is shown in
\cite{JacodShiryaev}, Theorem III 3.4.
The $n$-dimensional Poisson process $N$ has the martingale
representation property with respect to the filtration; see~\cite{JacodShiryaev}, Theorem III 4.37, so there exist
predictable processes $u^1, \ldots, u^n$ such that:
\[
S_t = 1 + \sum_{i = 1}^n \int_0 ^t u^i_s \,\mathrm{d}\overline{N}{}^i_s,
\]
where $\overline{N}{}^{i}_s := N^i_s - s$.

Now, let
\[
\mu^i_s := I ( S_{s-} > 0) \biggl( \frac{u^i_s}{ S_{s-}} + 1 \biggr)
\]
and note that
%
\begin{equation} \label{likelihood}
1 + \sum_{i = 1}^n \int_0 ^t S_{s-} ( \mu^i_s - 1) \,\mathrm{d}\overline{N}{}^{i}_s
= 1 + \int_0^tI( S_{s-} > 0) \,\mathrm{d}S_{s} = S_t,\qquad Q\mbox{-a.s.}
\end{equation}
The last equality follows from \cite{JacodShiryaev}, Lemma III 3.6.

Let $M^{(i)}_t := N^i_t - \int_0^t \mu^i_s \,\mathrm{d}s$ and note that
since $\Delta\overline{N}{}^{i}_s$
is bounded, \cite{JacodShiryaev}, Lemma III 3.14, says that the
quadratic (co)variation process $[\overline{N}{}^{i}, S]$ has locally integrable
variation, so its compensator $\langle\overline{N}{}^{i}, S \rangle$ is
well defined.
We can compute that
\[
\langle\overline{N}{}^{i} , S \rangle_t
= \int_0^t S_{s-} ( \mu_s^i - 1 )
\,\mathrm{d}\langle\overline{N}{}^{i} , \overline{N}{}^{i} \rangle_s =
\int_0^t S_{s-} ( \mu_s^i - 1 ) \,\mathrm{d}s,
\]
so we get from Girsanov's theorem (see \cite{JacodShiryaev}, Lemma III
3.14) that
\[
\overline N_t^i - \int_0 ^t \frac{ 1} { S_{s-}} \,\mathrm{d} \langle\overline{N}{}^{i} , S \rangle_s =
\overline{N}{}^{i}_t -
\int_0^t (\mu^i _s - 1)
\,\mathrm{d}s = M^{(i)}_t
\]
is a local $P$-martingale with respect to $\{\F_t\}_t$ for every $i
\leq n$.

Now
\[
M^{(i)}_t - M^i_t = \int_0 ^t \lambda_s - \mu^i _s \,\mathrm{d}s
\]
defines a continuous finite variation $P$-martingale, so $\mu^i =
\lambda^i$ $P$-a.s. a.e. and
$S = \mathcal{ E} ( K)$, where
$
K_t := \sum_{i = 1}^n \int_0 ^t ( \lambda_s^i - 1) \,\mathrm{d}\overline N_s^i
$
and
$\mathcal E$ is the stochastic exponential.
Let
$
K^C_t := \int_0 ^t (\lambda_s^C - 1 )\,\mathrm{d}\overline N^C_s
$
and
$
K^L_t := K^C_t - K_t.
$
Since $[K^L, K^C] = 0$ $Q$-a.s., we have that:
%
\begin{equation} \label{factorization}
\mathcal E (K
) = \mathcal E (K^C + K^L
)
= \mathcal E (K^C + K^L + [ K^C, K^L]
)
= \mathcal E (K^C
) \mathcal E (K^L
).
\end{equation}
The last equality follows from \cite{Protter}, Theorem II 38.

We now consider filtrations corresponding to the $\sigma$-algebras: $\G
_t := \F_{t \wedge T_D}^{A, D}$,
$\G_t^C := \F_{t \wedge T_D}^{A,C, D}$, $\G_t^L := \F_{t \wedge
T_D}^{A, D,L}$ and $\G_t^{C, L} := \F_{t \wedge T_D}$.
Moreover, we let $\tilde\lambda^D$ denote the
$\G_t$-predictable projection of $\lambda^D$; see \cite{JacodShiryaev}, Theorem
I 2.28. It is the unique $\G_t$-predictable process
such that\looseness=-1
\[
E [ \lambda_S^D | \G_{S -} ] =
\tilde\lambda_S^D
\]
for every $\G_t$-predictable stopping time $S$.

The local independence relations:
\begin{enumerate}[(1)]
\item[(1)]$L \nrightarrow_{T_D} C | A\cup C$;
\item[(2)]$C \nrightarrow_{T_D} L | A\cup D$;
\item[(3)]$C \nrightarrow_{T_D} A | L\cup D$;
\item[(4)]$C \nrightarrow_{T_D} D | A\cup L$
\end{enumerate}
and \eqref{factorization}
provide a factorization
\[
\frac{ \mathrm{d} P|_{ \G_{t-}^{C, L}}} {\mathrm{d} Q |_{ \G_{t-}^{C, L}} } =
S^L_t \cdot S^C_t,
\]
where $S^L$ is $\G_t^L$-predictable and
$S^C$ is $\G_t^C$-predictable. Bayes' theorem now gives that, whenever~$F$ is $\G_{t-}^L$-measurable and bounded, then
%
\begin{equation} \label{invariance}
E [ F | \G_{t-} ] = E [ F | \G_{t-}^C ]
\end{equation}
$P$-a.s.
Therefore, if we let $F$ be bounded and $\G_t^C$-predictable, then we
can compute:
\begin{eqnarray*}
E \biggl[ \int_0 ^T F_s \tilde\lambda_s^D \,\mathrm{d}s \biggr] & =& \int_0 ^T E
[ F_s \tilde\lambda_s^D ] \,\mathrm{d}s = \int_0 ^T E
[ F_s E [ \lambda_s^D| \G_{s-} ] ] \,\mathrm{d}s
\\
& =& \int_0 ^T E
[ F_s E [ \lambda_s^D| \G_{s-}^C ] ] \,\mathrm{d}s
= \int_0 ^T E
[ E [ F_s \lambda_s^D| \G_{s-}^C ] ] \,\mathrm{d}s
\\
& =& E \biggl[ \int_0 ^T F_s \lambda_s^D \,\mathrm{d}s \biggr].
\end{eqnarray*}
If we let $\tilde M_t^D = D_t - \int_0 ^t \tilde\lambda_s^D \,\mathrm{d}s$ and
$M_t^D = D_t -\int_0 ^t \lambda_s^D \,\mathrm{d}s $, then we
can compute:
\begin{eqnarray*}
E \biggl[\int_0 ^T F_s \,\mathrm{d} \tilde M_s^D\biggr] & = &
E \biggl[\int_0 ^T F_s \,\mathrm{d}D_s\biggr] - E \biggl[\int_0 ^T F_s \tilde\lambda_S^D \,\mathrm{d}s\biggr] \\
& =& E \biggl[\int_0 ^T F_s \,\mathrm{d}
D_s\biggr] - E \biggl[\int_0 ^T F_s \lambda_S ^D\,\mathrm{d}s\biggr] \\
& =& E \biggl[\int_0 ^T F_s \,\mathrm{d} M_s^D \biggr] = 0,
\end{eqnarray*}
so $\tilde M^D $ is a $P$-martingale with respect to the filtration
$\G_t^C$.
Now, since $\G_t \subset\G_t^L$, we have that $\tilde\lambda^D$ is $\G
_t^L$-predictable:
\begin{eqnarray*}
\tilde M^D_t & =& E [ \tilde M^D_{T_D} | \G_t^C ] = E
[ \tilde M^D_{T_D} I( T < t ) + \tilde M^D_{T_D} I( T \geq
t ) | \G_t^C ] \\
& =& \tilde M^D_{T_D} I( T < t ) + E
[ \tilde M^D_{T_D} I( T \geq
t ) | \G_t^C ] \\ &
=& \tilde M^D_{T_D} I( T < t ) + E
[ \tilde M^D_{T_D} | \F_t^{A, C, D} ] I( T \geq
t ) \\ & =& E [ \tilde M^D_{T_D} | \F_t^{A, C, D}
],
\end{eqnarray*}
that is, $\tilde M^D$ is a $P$-martingale with respect to the filtration
$\{\F_t^{A, C, D}\}_t$.

Finally, we note that
%
\begin{equation} \label{simple}
N^A_{t \wedge T_D} - \int_0 ^t \lambda_s^A I (s \leq T_D) \,\mathrm{d}s
\end{equation}
defines a $P$-martingale with respect to $\F_t$. Since $\lambda_s^A$
is $\F_t^{A,D}$-predictable, we also see that~\eqref{simple} defines a
martingale with
respect to $\{\F_t^{A, C, D}\}_t$.
\end{pf}

\section{Existence of randomized trial measures}

We have now come to the construction of randomized trial measures.
The idea is to
construct a reasonable randomized trial measure $\tilde P$ from the
observational measure $P$ such that $\tilde P \ll P$. The absolute
continuity is important since this provides a
natural method for simulating the empirical expectation of random
variables as if the data was sampled from the
counterfactual trial, while actually using $P$-distributed samples.
To get an idea of how this is done,
let $J \in\mathbb N$, let $H$ be a bounded random variable and let
$\omega_1,
\ldots, \omega_J$ be $J$ independently $P$-distributed samples from
$\Omega$. The law of large numbers then yields:
\[
\lim_{J \rightarrow\infty} \frac{ 1}{J} \sum_{j = 1}^J \frac{\mathrm{d}\tilde P
}{ \mathrm{d} P
}( \omega_j) H( \omega_j) = E_P \biggl[ \frac{\mathrm{d}\tilde P }{ \mathrm{d} P
} H \biggr] = E_{\tilde P} [ H ], \qquad P  \mbox{-a.s.}
\]
Heuristically, this means that the likelihood ratio can be viewed as
a transformation from the observational
study into the counterfactual scenario.

There might exist several reasonable counterfactual trials, each
corresponding to a~choice of a well-behaved
treatment and censoring strategies. Given a non-negative
$\F_t^{A,D}$-predictable process~$\tilde\lambda^A$ and a non-negative
$\F_t^{A,C,D}$-predictable process $\tilde\lambda^C$, we can consider the
problem of finding a randomized trial measure $\tilde P$ that has
$\tilde\lambda^A$ as the $\F_t$-intensity of $N^A$ and $\tilde
\lambda^C$ as the $\F_t$-intensity of $N^C$. This suggests that
$\tilde\lambda^A$ is the treatment strategy and $\tilde\lambda^C$
is the censoring strategy in the counterfactual trial.
We will consider the counterfactual treatment strategy given by the
$P$-intensity of $N^A$ with respect to $\F_t^{A,D}$.
The counterfactual censoring strategy will be given by the
$P$-intensity of $N^C$ with respect to $\F_t^{A,D,C}$.
This gives a randomized trial measure with a likelihood
ratio that heuristically corresponds to the \textit{stabilized weights} one
usually considers in the discrete time marginal structural models; see
\cite{Robins1}.
The problem of finding such a randomized trial measure is a martingale problem.
Note that this problem might not have a solution.
The next theorem shows that if the counterfactual
strategies are not too different from the observed
intensities, then there exists a unique corresponding randomized trial
measure~$\tilde P$.

\begin{theorem} \label{maintheorem}
Suppose that there exist positive numbers $\theta_1$ and $\theta_2$
such that:
%
\begin{equation} \label{mainconditionA}
\lambda^A_s - \theta_1 \sqrt{ \lambda^A_s }
\leq
E_P[\lambda^A_s
|\F_{s-}^{A, D} ] \leq \lambda^A_s + \theta_1 \sqrt{
\lambda^A_s }
\end{equation}
and
%
\begin{equation} \label{mainconditionC}
\lambda^C_s - \theta_2 \sqrt{ \lambda^C_s }
\leq
E_P[\lambda^C_s
|\F_{s-}^{A,C, D} ] \leq \lambda^C_s + \theta_2 \sqrt{
\lambda^C_s }
\end{equation}
for almost every $s$ $P$-a.s. Let $\tilde\lambda^A$ denote the
$P$-intensity of $N^A$ with respect
to the filtration $\{\F_t^{A, D}\}_t$ and let
$\tilde\lambda^C$ denote the $P$-intensity of $N^C$ with respect
to the filtration $\{\F_t^{A,C, D}\}_t$.

The equation
\begin{eqnarray}
\label{likelihoodratio}
R_t &:=& \prod_{ s \leq t} \biggl( \frac{ \tilde\lambda_s^A }{
\lambda_s^A } \biggr)^{\Delta N^A_s }
\exp\biggl( \int_0^t \tilde \lambda_s^A - \lambda_s^A \,\mathrm{d}s
\biggr)
\nonumber
\\[-8pt]
\\[-8pt]
\nonumber
&&{}\times
\prod_{ s \leq t} \biggl( \frac{ \tilde\lambda_s^C }{
\lambda_s^C } \biggr)^{\Delta N^C_s }
\exp\biggl( \int_0^t \tilde \lambda_s^C - \lambda_s^C \,\mathrm{d}s \biggr)
\end{eqnarray}
defines a square-integrable $P$-martingale with respect to the
filtration $\{\F_t\}_t$.
Moreover,
\[
\mathrm{d} \tilde P = R_T \,\mathrm{d}P
\]
defines a randomized trial measure on $(\Omega, \F)$ such that
the martingale dynamics of the biological processes $D$ and $L$
coincide for the two probability measures $\tilde P$
and $P$, that is,
\[
L_t - \int_0 ^t H_s^L \lambda_s^L \,\mathrm{d}s   \quad \mbox{and}\quad   D_t - \int_0 ^t
\lambda_s^D I( s \leq T_D) \,\mathrm{d}s
\]
also define $\tilde P$-martingales with respect to $\{\F_t\}_t$.
\end{theorem}

\begin{pf}
We define
\[
K_t := \int_0 ^t \biggl(\biggl( \frac{ \tilde\lambda_s ^A} { \lambda_s ^A} - 1
\biggr) \,\mathrm{d}M^A_s
+
\int_0 ^t \biggl( \frac{ \tilde\lambda_s ^C} { \lambda_s ^C} - 1
\biggr) \biggr)\,\mathrm{d}M^C_s.
\]
By the innovation theorem, we have that
%
\begin{equation}\label{innovation}
\tilde\lambda_s^A = E_P[\lambda^A_s
|\F_{s-}^{A, D} ]  \quad \mbox{and}\quad
\tilde\lambda_s^C = E_P[\lambda^C_s
|\F_{s-}^{A, C,D} ] \qquad P  \mbox{-a.s., } s \mbox{-a.e.}
\end{equation}
By \eqref{mainconditionA} and
\eqref{mainconditionC}
we have that
\[
\biggl( \frac{ \tilde\lambda_s ^A} { \lambda_s ^A} - 1
\biggr)^2 \lambda_s^A +
\biggl( \frac{ \tilde\lambda_s ^C} { \lambda_s ^C} - 1
\biggr)^2 \lambda_s^C
\leq\theta_1 + \theta_2,\qquad  P \mbox{-a.s., }   s \mbox{-a.e.}
\]
We therefore obtain that:
\[
\langle K, K \rangle_t  = \int_0 ^t \biggl(\biggl( \frac{ \tilde\lambda_s ^A}
{ \lambda_s ^A} - 1
\biggr)^2 \lambda_s^A + \biggl( \frac{ \tilde\lambda_s ^C} { \lambda_s
^C} - 1
\biggr)^2 \lambda_s^C\biggr) \,\mathrm{d}s \leq
( \theta_1 + \theta_2 ) \cdot t.
\]

Since $\langle K, K \rangle$ is bounded on the interval $[0, T]$, \cite
{Lepingle}, Theorem
II.1, yields that
the stochastic exponential $R := \mathcal E ( K)$, given by
the SDE:
\begin{equation}
\label{Doleansexp}
R_t = 1 + \int_0 ^t R_{s-} \,\mathrm{d}K_s
\end{equation}
is square-integrable.

Now
\[
\tilde \mathrm{d}P = R_T \,\mathrm{d}P
\]
defines a probability measure $\tilde P$ on $(\Omega,\F)$.

Note that we have:
\[
M^D - \int_0^\cdot \frac{1 }{R_{s-} } \,\mathrm{d}\langle M^D, R \rangle_s = M^D
\quad \mbox{and}\quad
M^l - \int_0^\cdot \frac{1 }{R_{s-} }\, \mathrm{d}\langle M^l, R \rangle_s = M^l
\]
$P$-a.s. for every $l \in L$. Moreover, Girsanov's theorem (see \cite
{JacodShiryaev}, Lemma III 3.14)
gives that these processes define local $\tilde P$-martingales with
respect to the filtration $\{\F_t\}_t$.

Moreover,
\[
\langle M^A , R \rangle_t
= \int_0^t R_{s-} \biggl( \frac{ \tilde\lambda_s^A
}{\lambda_s^A } - 1 \biggr) \,\mathrm{d}\langle M^A, M^A\rangle_s =
\int_0^t R_{s-} ( \tilde\lambda_s ^A - \lambda_s ^ A
) \,\mathrm{d}s,
\]
so again by Girsanov's theorem, we have that
\[
M_t^A - \int_0 ^t \frac{ 1} { R_{s-}} \,\mathrm{d} \langle M^A, R \rangle_s =
M_t^A -
\int_0^t \tilde\lambda_s ^A - \lambda_s ^ A
\,\mathrm{d}s =
N^A_t - \int_0 ^1\tilde \lambda_s^A \,\mathrm{d}s
\]
defines a $\tilde P$-martingale with respect to the filtration $\{\F_
t\}_t$.
Analogously, we have that
\[
M_t^D - \int_0 ^t \frac{ 1} { R_{s-}} \,\mathrm{d} \langle M^D, R \rangle_s
=
N^D_t - \int_0 ^t\tilde \lambda_s^D \,\mathrm{d}s
\]
defines a local $\tilde P$-martingale with respect to the filtration $\{
\F_
t\}_t$.\vadjust{\goodbreak}

Finally, we note that by \cite{JacodShiryaev},
Theorem I 4.60, the SDE \eqref{Doleansexp} has the explicit solution
given by~\eqref{likelihoodratio}. Expressions of this form are well
known in the literature on marked point process; see, for
instance, \cite{Jacod}.
\end{pf}
\begin{remark}
Note that the condition \eqref{mainconditionA} heuristically means that the
short-term ``risk'' of starting treatment, given the previous full
history,
\[
\lim_{h \rightarrow0} h^{-1} P ( t \leq T_A < t + h |
\F_t ),
\]
is not too different from the
short-term ``risk'' of starting treatment when we do not pay
attention to the underlying health process or censoring, that is,
\[
\lim_{h \rightarrow0} h^{-1} P ( t \leq T_A < t + h |
\F_t^{A, D} ).
\]

Similarly, condition \eqref{mainconditionC} heuristically means that the
short-term ``risk'' of being censored, given the previous full
history,
is not too different from the
short-term ``risk'' of being censored when we do not pay
attention to the underlying health process.
In words, this means that
the previous history of the health process $L$ alone can not at any
time yield too high of a short-term
``risk'' of starting treatment or being censored.

\end{remark}

\section{Weighted additive hazard regression}

Suppose that we have observations of
$m$ independent individuals until death or censoring
from an observational study and want to estimate the total effect of
treatment.
Ideally, we would like to base our estimate on some randomized
trial. However, such a trial might not be available. The marginal structural
approach now suggests that we simulate a counterfactual randomized trial,
using the data we already have.
We would then like to estimate the counterfactual hazard,
that is, the hazard that the patient would have if he, contrary to the fact,
had participated in the
randomized trial. In this way we could compare the
total effect of being in treatment versus never being in treatment.

We assume that the observations from each individual are
$P$-distributed. The observations of the patient would have been
$\tilde P$-distributed if he, contrary to the fact, participated in the
hypothetical randomized trial.

We assume that the counterfactual intensity follows Allen's additive hazard
regression model, see \cite{Andersen,AalenGjessingBorgan}.
To formalize this, let $\beta^0$ and $\beta^1$ be functions on $[0,
T]$.
We assume that we only have instantaneous effect and that the hazard
for the event, with
respect to the treatment and event history, is given by:
\[
\beta_t^0 + \beta^1_t A_{t-}.
\]
One way to estimate the hazard from the counterfactual trial is to
weight the observations
by the corresponding likelihood ratios.
We will prove that a suitable weighted
variant of Aalen's additive hazard regression gives
consistent estimators of $\int_0^t \beta^0_s \,\mathrm{d}s$ and $\int_0^t
\beta^1_s \,\mathrm{d}s$ from independent $P$-distributed observations.
This requires some notation. Let $Y_s^1, \ldots, Y^m_s$ be the
``at-risk'' indicators and
$A_{t}^1, \ldots, A^m_{t}$ be the ``at-treatment'' indicators for the
$m$ independent individuals.
We define the $m \times2$-matrix:
\[
X^{(m)}_t =
\pmatrix{
Y^{1}_{t} & Y_{t}^{1} \cdot A_{t-}^1 \vspace*{2pt}\cr
\vdots& \vdots\vspace*{2pt}\cr
Y_t^m & Y_t^m \cdot A_{t-}^m
}
.
\]
Moreover, let
$R^1_t, \ldots, R^m_t$ be the individual likelihood ratios at time $t$
and let
\[
R_{t-} \m=
\pmatrix{
Y^1_t R^1_{t-} & 0 & \ldots& 0 \vspace*{2pt}\cr
0 & Y_t^2 R^2_{t-} & & \vdots\vspace*{2pt}\cr
\vdots& & \ddots& 0\vspace*{2pt} \cr
0 & 0 & \ldots& Y^m_t R^m_{t-}
}
.
\]
Finally, let $D_t^1, \ldots, D^m_t$ be event processes for the
$m$ individuals before $t$.
The observed events are now given by the vector
\[
D^{(m)}_t=
\pmatrix{
\displaystyle\int_0 ^t Y_s^1 \,\mathrm{d}D^1_s\vspace*{2pt}\cr
\vdots\vspace*{2pt}\cr
\displaystyle\int_0 ^t Y_s^m\, \mathrm{d}D^m_s
}
.
\]
\begin{theorem}
We assume that:
\begin{enumerate}[(1)]
\item[(1)] \label{condbound}
The $P$-intensity of $D$ with respect to the filtration
$\F_t$,
$\lambda^D$ is dominated
by an integrable function $G$.
\item[(2)](Positivity) Both the ``at-risk'' groups in the
counterfactual trial are always present, that is,
\[
E_{\tilde P} [
Y_sA_{s-} ] > 0 \quad \mbox{and}\quad  E_{\tilde P} [ Y_s(1 -
A_{s-}) ] > 0
\]
for every $s \in[0, T]$.
\item[(3)]
There exist integrable and left-continuous functions with right limits
$\beta^0$ and $\beta^1$ such that $\beta= ( \beta^0 , \beta^1)^{\mathrm{T}}$ and
such that
\[
D_t\m- \int_0 ^t X_s \m\beta_s \,\mathrm{d}s
\]
is a $\tilde P$-martingale with respect to the filtration
$\F_t^{A, C, D}$, i.e.,
$Y_t ( \beta_t^0 + \beta_t^1 A_{t-}^1)$ is the \mbox{$\tilde P$-intensity}
of $ D$ w.r.t. the filtration $\F_t^{A, C, D}$.
\end{enumerate}
We let
\begin{eqnarray*}
J_s^{(m)} & := &I \Biggl( \sum_{i = 1} ^m R^i_{t-} Y^i_t ( 1- A^i_{t-}) >
0 \mbox{ and } \sum_{i= 1}^m R^i_{t-} Y^i_t A^i_{t-} > 0\Biggr), \\
\hat B_t^{(m)} & := &\int_0^t J_s \m\bigl( X_s^{(m)\mathrm{T}} R_{s-}\m X_s\m\bigr) ^{-1}X_s
^{(m)\mathrm{T}} R\m_{s-}
\,\mathrm{d} D_s\m
\end{eqnarray*}
and
\[
B_t := \int_0 ^t \beta_s \,\mathrm{d}s.
\]

Now $ \hat B\m$ is a consistent estimator of $B_t$
in the sense that:
\[
\lim_m P\bigl( d \bigl( \hat B\m, B \bigr) \geq\epsilon\bigr) = 0
\]
for every $\epsilon> 0$, where $d$ denotes the Skorokhod metric; see
\cite{JacodShiryaev} or \cite{Billingsley}.
\end{theorem}

\begin{pf}
We define:
\[
Y^{(m)}_t=
\pmatrix{
Y^1_t\vspace*{2pt}\cr
\vdots\vspace*{2pt}\cr
Y^m_t
}
,\qquad
\lambda^{(m)}_t=
\pmatrix{
\lambda^1_t\vspace*{2pt}\cr
\vdots\vspace*{2pt}\cr
\lambda^m_t
}
\quad
\mbox{and}\quad
M^{(m)}_t=
\pmatrix{
D^1_t - \displaystyle\int_0 ^t \lambda_s^1 \,\mathrm{d}s\vspace*{2pt}\cr
\vdots\vspace*{2pt}\cr
D^m_t - \displaystyle\int_0 ^t \lambda_s^m \,\mathrm{d}s
}
.
\]
We will often drop the index $(m)$ in order to simplify the
notation. Another simplification of the notation we will use is
$E[\cdot ]$ for the expectation with
respect to $P$ and $\tilde E[\cdot]$ for the expectation with
respect to $\tilde P$.

First we prove that
%
\begin{equation} \label{large}
\lim_m P\biggl( \sup_{t \leq T} \biggl|
\int_0 ^t J_s(X_s^{\mathrm{T}} R_{s-} X_s)^{-1} X^{\mathrm{T}}_s R_{s-}
\lambda_s \,\mathrm{d}s- B_t \biggr| \geq\epsilon\biggr) = 0
\end{equation}
for every $\epsilon> 0$.
Define:
\[
V :=
\pmatrix{
1 & 0 \vspace*{2pt}\cr
-1 & 1
}
\quad \mbox{and}\quad
S_t :=
\pmatrix{
\displaystyle\sum_{i = 1}^m R^i_{t-} Y^i_t ( 1- A^i_{t-}) & 0
\vspace*{2pt}\cr
0 & \displaystyle\sum_{i = 1} ^m R^i_{t-} Y^i_t A^i_t
}
.
\]
The matrix $V$ is invertible and, using the fact that the $Y^i$ and the
$A^i$ are indicators, we have:
\[
V^{\mathrm{T}} X^{\mathrm{T}}_t R_{t-} X_t V = S_t,
\]
that is, $X^{\mathrm{T}}_t R_{t-} X_t$ is congruent to the diagonal matrix $S_t$.
A simple matrix computation gives that
\[
(X^{\mathrm{T}}_t R_{t-} X_t)^{-1} = V S^{-1}_t V^{\mathrm{T}},
\]
when $J_s > 0$.

Now, we see that:
\[
J_t (X_t^{\mathrm{T}} R_{t-} X_t)^{-1} X^{\mathrm{T}}_t R_{t-}
\lambda_t
= J_t V
\begin{pmatrix}
H_t^0 \\
H_t^1
\end{pmatrix}
,
\]
where

\[
H_t^0 = \frac{\sum_{i = 1} ^m R^i_{t-} Y_t^i ( 1 - A^i_{t-})
\lambda_t^i }{ \sum_{i = 1} ^m R^i_{t-} Y_t^i ( 1 - A^i_{t-}) }
\quad \mbox{and}\quad
H_t^1 = \frac{\sum_{ i = 1} ^m R^i_{t-} Y_t^i A^i_{t-}
\lambda_t^i }{ \sum_{ i= 1} ^m
R^i_{t-} Y_t^i A^i_{t-} } .
\]

Since $\tilde E[Y_tA_{t-}] > 0$, the law of large numbers implies that
$H_t^0$ converges in probability to
\[
\frac{\tilde E[ Y_t ( 1 - A_{t-}) \lambda_t^D] }{ \tilde E[
Y_t ( 1 - A_{t-}) ] } =
\tilde E[ \lambda_t^D | Y_t ( 1 - A_{t-} ) = 1].
\]
Analogously, since $E[Y_t(1-A_{t-})] > 0$, we have that
$H^1_t$ converges in probability to
\[
\frac{\tilde E[ Y_t A_{t-} \lambda_t^D ] }{
\tilde E[
Y_t A_{t-} ] } = \tilde E[ \lambda_t^D | Y_t A_{t-1} = 1].
\]
By a similar argument, we see that $\{J_s^{(m)} \}$ converges
in probability to $1$ for almost every~$s$.

Since the $P$-intensity of $D$ with respect to $\F_s^{A, C ,D}$
coincides with $E[ \lambda^D_s | \F_{s-}^{A, C, D}]$ $P$-a.s. for
almost every $s$, we have that:
\[
\tilde E[ \lambda_D | Y_t ( 1 - A_{t-} ) = 1] = \beta_t^0
\quad \mbox{and}\quad
\tilde E[ \lambda_t^D | Y_t A_{t-} = 1] = \beta_t^0 + \beta_t^1.
\]
This means that
\begin{equation} \label{integrand}
J_s (X_t^{\mathrm{T}} R_{t-} X_t)^{-1} X^{\mathrm{T}}_t R_{t-} \lambda_t
\end{equation}
converges in probability to
$\beta_t$ when $m$ increases. Note that
\[
E \Biggl[\sup_t\biggl | \int_0 ^t J_s V
\pmatrix{ H_s^0 \vspace*{2pt}\cr
H_s^1
}
\,\mathrm{d}s - \int_0 ^t \beta_s \,\mathrm{d}s \biggr|
\Biggr] \leq
\int_0 ^{T} E\Biggl[ \biggl| J_s V
\pmatrix{ H_s^0 \vspace*{2pt}\cr
H_s^1
}
- \beta_s \biggr| \Biggr] \,\mathrm{d}s,
\]
so by the dominated convergence theorem, we obtain \eqref{large}.

We will now prove that
\begin{eqnarray*}
Z_t & :=& \hat B_t -
\int_0 ^t J_s (X_s^{\mathrm{T}} R_{s-} X_s)^{-1} X^{\mathrm{T}}_s R_{s-}
\lambda_s \,\mathrm{d}s
\\
&= &
\int_0 ^t (X_s^{\mathrm{T}} R_{s-} X_s)^{-1} X^{\mathrm{T}}_sR_{s-} \,\mathrm{d}M_s
\end{eqnarray*}
converges weakly to~$0$.
Note that:
\begin{eqnarray*}
\langle Z , Z \rangle_t
& =& \int_0 ^t J_s ( X_s^{\mathrm{T}} R_{s-} X_s ) ^{-1} X_s ^{\mathrm{T}} R_{s-} \,\mathrm{d} \langle
M, M\rangle_s R_{s-} X_s
( X_s^{\mathrm{T}} R_{s-} X_s ) ^{-1} \\
& =& \int_0 ^t J_s V S_s^{-1} V^{\mathrm{T}} X_s ^{\mathrm{T}} R_{s-} \,\mathrm{d} \langle M, M \rangle_s
R_{s-} X_s V S_s^{-1} V^{\mathrm{T}}
\\
& = &V \int_0 ^t J_s S_s ^{-2} U_s\,\mathrm{d}s V^{\mathrm{T}},
\end{eqnarray*}
where
\[
U_s =
\pmatrix{ \displaystyle\sum_{i = 1} ^m R^{i 2}_{s-} Y_s^i
( 1 - A^i _{s-} ) \lambda_s^i & 0 \vspace*{2pt}\cr
0 & \displaystyle\sum_{i = 1} ^ n R^{i 2}_{s-} Y_s^i A^i _{s-} \lambda_s^i
}
.
\]
Now,
\[
E \Biggl[ \sup_t \biggl| \int_0^t V J_s S^{-2}_s U_s V^{\mathrm{T}} \,\mathrm{d}s \biggr| \Biggr]
\leq\int_0 ^{T} E \bigl[ | V J_s S_s ^{-2} U_s V^{\mathrm{T}} | \bigr] \,\mathrm{d}s,
\]
so, by the dominated convergence theorem,
$\{\langle Z\m, Z\m\rangle\}$ converges uniformly in probability to~$0$.

We define
\[
Z^{(\epsilon,m)}_t := \int_0 ^t I \bigl(
| J_s (X_s^{\mathrm{T}} R_{s-} X_s ) ^{-1} X^{\mathrm{T}} R_{s-} Y_s
| \geq \epsilon
\bigr) J_s
(X_s^{\mathrm{T}} R_{s-} X_s ) ^{-1} X^{\mathrm{T}}R_{s-} \,\mathrm{d}M_s
\]
and see that
\[
0 \leq\bigl\langle Z^{(\epsilon,m)}, Z^{(\epsilon,m)} \bigr\rangle_t \leq
\bigl\langle Z\m, Z\m
\bigr\rangle_t.
\]
Since both
\[
\bigl\{ \bigl\langle Z^{(\epsilon,m)}, Z^{(\epsilon,m)} \bigr\rangle\bigr\}_m \quad \mbox{and}\quad  \bigl\{\bigl\langle Z\m, Z\m
\bigr\rangle\bigr\}_m
\]
converge uniformly in probability to $0$, the central limit
theorem for martingales (\cite{JacodShiryaev}, Theorem~VIII 3.22)
implies that $\{Z^{(m)}\}_m $
converges weakly to $0$.

We have that
\[
\hat B_t ^{(m)} = Z_t^{(m)}
+
\int^t_0 J_s (X_s^{\mathrm{T}} R_{s-} X_s)^{-1} X^{\mathrm{T}}_s R_{s-} \lambda_s\m
\,\mathrm{d}s,
\]
so the sequence
$\{ \hat B^{(m)} \}_m $
is the sum of two $C$-tight sequences.
Jacod and Shiryaev \cite{JacodShiryaev}, Corollary VI
3.33, implies that the sequence itself is also $C$-tight. By
Slutsky's theorem, the finite-dimensional distributions\vspace*{-2pt} on the form
$\mathcal L ( \hat B^{(m)}_{t_1} , \ldots, \hat B^{(m)}_{t_j} )$
converge weakly to the Dirac measures:
$
\delta_{B_{t_1}, \ldots, B_{t_j} }.
$
This means that $\{ \hat B^{(m)} \}_m $ converges in law and
therefore in probability to $B$.
\end{pf}
%

\section{Concluding remarks}
We have shown that marginal structural modeling can be understood in
terms of change of probability measures. The author believes that this
is an elucidating point of view that is natural in the
framework of modern probability theory.

As stressed by several authors, there is a very important and highly non-trivial
assumption one has to make in order to interpret effects from marginal
structural
models as causal.
This is the assumption of \textit{no unmeasured confounders}, or
equivalently: \textit{all confounders are measured}.
This means that every process that affects the short-term behavior of
both the treatment and the censoring or both the treatment and the
event must be observed.
In this equivalent form, it becomes more apparent that this is just an
assumption about
completeness of the model. This completeness assumption is not that mysterious.
When modeling various phenomena in the natural sciences, one typically
assumes that all the important variables are contained in the
model. This is also necessary in the MSM approach.
However, it is important to note that this is not generally a
statistically testable assumption.
It is also not a condition that would follow from a~%
mathematical argument without further assumptions about the model.

Heuristically, the MSM
approach provides an adjustment of the the treatment
effect bias caused by the measured confounders.
In the marginal structural model approach, instead of
modeling the underlying and potentially very complicated biology,
one models a randomized trial. The problem of computing
marginal effects then splits into two parts.
The first problem is to model
the marginal intensity of the event in the simulated ``randomized
trial''.
If one knew the corresponding likelihood ratio process, then
this would be obtainable using, for instance, the weighted additive hazard
regression from the previous section.
In order to compute this likelihood ratio, one has to
deal with the second part of our problem. That is to model
the dynamics of the treatment and censoring processes
given the full and the marginal history in the observational
study.
This is a crucial point. We have chosen not to deal with this
problem in the current paper. However, one could use regression
techniques to do this at least approximately.
In the discrete time setting, one typically uses pooled logistic
regressions see \cite{Sterne,Zidovudine}. In the
continuous-time setting it is probably more natural to use additive
hazard or Poisson regression to estimate the censoring and treatment
intensities, both with respect to the full covariate history
and marginal covariate history. This will be the topic of future
work. Once these intensities are known, one can compute the likelihood
ratio process using \eqref{likelihoodratio}.

\section*{Acknowledgements}
Supported by the Research Council of Norway. Project: 170620/V30.
I would like to thank Odd~O. Aalen, Vanessa Didelez and Jon Michael
Gran for helpful discussions related to this project.

\printhistory

\end{document}